\documentclass[twoside,12pt]{amsart}
\usepackage{amsmath,latexsym,amssymb,times,mathptm,amsfonts,enumerate}

\def\proof{\noindent{\bf{Proof.} }}
\def\sqr#1#2{{\vcenter{\hrule height.#2pt
        \hbox{\vrule width.#2pt height#1pt \kern#1pt
                \vrule width.#2pt}
        \hrule height.#2pt}}}
\def\square{\mathchoice\sqr64\sqr64\sqr{4}3\sqr{3}3}
\def\QED{\hfill$\square$\\}

\def\tratto{\mbox{\rule{2mm}{.2mm}$\;\!$}}

\newtheorem{theorem}{Theorem}[section]

\newtheorem{corollary}[theorem]{Corollary}
\newtheorem{lemma}[theorem]{Lemma}
\newtheorem{proposition}[theorem]{Proposition}

\newtheorem{remark}[theorem]{Remark}

\newtheorem{example}[theorem]{Example}

\newtheorem{assumptions}[theorem]{Assumptions and Discussion}
\newtheorem{notation}[theorem]{Notation}

\newcommand{\m}{\mathfrak{m}}
\newcommand{\f}[1]{\ensuremath{\mathfrak{#1}}}
\newcommand{\ol}[1]{\ensuremath{\overline{#1}}}

\newcommand{\depth}{\ensuremath{{\rm{depth}}\;}}
\newcommand{\core}[1]{\ensuremath{{\rm{core}}(#1)}}

\numberwithin{equation}{section}

\begin{document}
\baselineskip=16pt

\title{Minimal reductions and cores of edge ideals}
\author[L. Fouli]{Louiza Fouli}
\address{Department of Mathematical Sciences,  New Mexico State University,  Las Cruces, New Mexico 88003, USA}
\email{lfouli@math.nmsu.edu}

\author[S. Morey]{Susan Morey}
\address{Department of Mathematics, Texas State University, 601 University Drive, San Marcos, Texas 78666, USA}
\email{morey@txstate.edu}

\subjclass[2010]{13A30, 13A15, 05E40}
\keywords{edge ideal, even cycles, minimal reductions, core}

\begin{abstract} 
We study minimal reductions of edge ideals of graphs and determine  restrictions on the coefficients of the generators of these minimal reductions. We prove that when $I$ is not basic, then $\core{I}\subset \m I$, where $I$ is an edge ideal in the corresponding localized polynomial ring and $\m$ is the maximal ideal of this ring. We show that the inclusion is an equality for the edge ideal of an even cycle with an arbitrary number of whiskers. Moreover, we show that the core is obtained as a finite
intersection of homogeneous minimal reductions in the case of even cycles. The formula for the core does not hold in general for the edge ideal of any graph and we provide a counterexample. In particular, we show in this example that the core is not obtained as a finite intersection of general minimal reductions.

\end{abstract}

\maketitle

\section{Introduction}\label{intro}

Let $R$ be a Noetherian ring and $I$ an ideal of $R$. Recall that a
{\it reduction} 
of $I$ is an ideal $J$ such that $J \subset I$ and $\ol{I}=\ol{J}$,
where $^{^{\tratto}}$ denotes the integral closure. Equivalently, $J
\subset I$ is a reduction of $I$ if and only if $I^{r+1}=JI^{r}$ for some 
nonnegative integer $r$ \cite{NR}.  When $R$ is a
Noetherian local ring then we may consider minimal reductions, where
minimality is with respect to inclusion. Northcott and Rees
proved that 
when $R$ is a Noetherian local ring with infinite residue field then either
$I$ has infinitely many minimal reductions or $I$ is {\it basic},
i.e. $I$ is the only reduction of itself. 

A reduction can be thought of as a simplification of the ideal. One advantage to considering reductions is that they are in principle smaller ideals with the same asymptotic behavior as the ideal $I$ itself. For example, all  minimal reductions  of $I$ have the same height, the same radical, and the same multiplicity as $I$. 

Let $R$ be a Noetherian local ring with infinite residue field and $I$ an ideal of $R$. Then every minimal reduction $J$ of $I$ has the same minimal number of generators, $\ell(I)$, where $\ell(I)$ is the {\it analytic spread} of $I$  (see Section~\ref{background}). It is well known that every minimal generating set of a reduction $J$ of $I$ can be extended to a minimal generating set of $I$. Therefore $\ell(I) \leq \mu(I)$, where $\mu(I)$ denotes the minimal number of generators of $I$. When $\ell(I)=\mu(I)$ then $I$ is basic.

Minimal reductions are not unique and therefore one considers the intersection of all the reductions of an ideal, namely the {\it core} of the ideal. This object was defined by Rees and Sally  \cite{RS}.  When $R$ is a Noetherian local ring
it is enough to consider the intersection of the minimal 
reductions. This intersection is in general infinite and there is significant difficulty in obtaining closed formulas that describe the core. Several authors have determined 
formulas that compute the core under various assumptions; Corso,  Huneke, Hury,  Polini, Smith,
 Swanson, Trung, Ulrich, Vitulli to name a few, \cite{CPU01, CPU02,HS1, HT, HySm1, HySm2,PU, PUV}.
 Furthermore, Hyry and Smith have discovered a connection with a
celebrated conjecture by Kawamata on the non-vanishing of sections
of line bundles \cite{HySm1}. They prove that the validity of the conjecture is equivalent to a statement about {\it gradedcore} and thus renewed the interest in understanding the core. The graded core is the intersection of all homogeneous minimal reductions and in general, ${\rm graded}\core{I} \subset \core{I}$. In Section~\ref{finiteintersection} we provide an instance where equality holds.

 In \cite{PUV} Polini, Ulrich and Vitulli
study the core of $0$-dimensional
monomial ideals in polynomial rings. They prove that  the core is obtained by computing the mono of a general locally minimal reduction of $I$ \cite[Theorem~3.6]{PUV}. The mono of an ideal $K$ is the largest monomial subideal contained in $K$. They provide an effective algorithm for computing the core, which is implemented in computer algebra programs such as CoCoA. In general, though, the question of what is the core of a monomial ideal is quite open.

 It was shown in \cite[Proposition~2.1]{Singla} that among the monomial reductions of a
monomial ideal, there is a unique minimal element. However, this
reduction need not be minimal among all reductions. If
  the monomial 
ideal $I$ has a square-free generating set, then Singla showed
that the only monomial reduction 
of $I$ is $I$ itself \cite[Remark~2.4]{Singla}. This leaves a large class of monomial ideals
whose minimal reductions are not monomial. Even though a monomial ideal need not have monomial minimal reductions its core is monomial  \cite[Remark~5.1]{CPU01}.

The class of square-free monomial ideals generated in degree two can be viewed as
edge ideals of graphs (see Section~\ref{background}).  Such ideals
were introduced in \cite{Rafael} and their  
properties have been studied by many 
authors, including \cite{AJ, dochtermann, FHVodd, katzman1, kummini, SVV, Vill}. In
order to discuss 
minimal reductions, the ring needs to be a local ring with infinite
residue field. Since $I$ is a homogeneous ideal, we will view $I$ as
an 
ideal in the localization of a polynomial ring at its homogeneous maximal ideal $\m$ and we will assume that the residue field is infinite. By abuse of
notation we will still 
denote the ideal by $I=I(G)$, where $G$ is the associated graph.   We note here that the edge ideals we study are far from being $0$-dimensional, so the monomial ideals we consider are not in the same class as the ones considered by Polini, Ulrich, and Vitulli in \cite{PUV}.

As mentioned earlier, $\ell(I) \leq \mu(I)$ and when $\ell(I)=\mu(I)$ then the ideal is basic. In this case the core is trivial, i.e. $\core{I}=I$.  When  $I$  is an ideal with  $\ell(I)=\mu(I)-1$ then $I$ is called an ideal of {\it second analytic deviation one}. For these ideals we show that if $(h_1, \ldots, h_s)$ is a minimal
generating set of $I$, then $J$ has a generating set of
the form   
$(h_1+a_1h_{t}, h_2+a_2h_{t}, \ldots, h_{s}+a_{s}h_{t})$ 
for some $1 \leq t \leq s$, where $a_i \in R$ for all $i$ and $a_t=-1$
(Lemma~\ref{gen red format}). In Corollary~\ref{gen format} we extend this to give a description of the structure of minimal reductions of any ideal in a Noetherian local ring. 
Not all choices of $a_i$ will result in a reduction, even when the second analytic deviation is one. One of the goals of this paper is to find restrictions on the coefficients $a_i$.  When $I$ is the edge ideal of a graph with a unique even cycle of length $d$ then $I$ is an ideal of second analytic deviation one (Remark~\ref{basic properties}). We show that if $\prod \limits_{i=1}^{{\frac{d}{2}}} a_{2i-1} = \prod   \limits_{j=1}^{{\frac{d}{2}}} a_{2j}$ then $J$ is not a reduction of $I$ (Corollary~\ref{productofbs}). The condition that $J$ is a minimal reduction of $I$ is  an open condition, i.e. the vectors of the coefficients $a_i$ are in a dense open subset of $\mathbb{A}_{R}^{s-1}$. More precisely, we show that there exists a hypersurface defined by the relation on the products of the coefficients $a_i$ as above, in the complement of this open set.

Let $I$ be the edge ideal of a graph that is not basic and let $R$ be the corresponding localized polynomial ring. Let $\m$ be the maximal ideal of $R$. We show in Theorem~\ref{core in mI} that $\core{I}\subset \m I$.
To establish a case where equality occurs, we consider the class of edge ideals of even cycles  with an arbitrary number of whiskers (potentially none) at each vertex. Let $I$ be such an ideal. We show that   $J:I=\m$ for all minimal reductions $J$ of $I$, Theorem~\ref{colonm full whiskers}. In particular, these results imply that $J:I$ is independent of the choice of the minimal reduction $J$ of $I$. This means that $I$ is a balanced ideal in the sense of \cite{U2}. This balanced property allows us to compute a formula for the core of these ideals.

Let $R$ be a Gorenstein local ring and let $I$ be an ideal of $R$ that satisfies $G_{\ell}$ and is weakly $(\ell -1)$-residually $S_2$, where $\ell=\ell(I)$. Under these assumptions Corso, Polini and Ulrich prove that $\core{I}=(J:I)J=(J:I)I$ for any minimal reduction $J$ of $I$ \cite[Theorem~2.6]{CPU02}. The edge ideals we consider are not weakly $(\ell-1)$-residually $S_2$. Nonetheless, we establish the same formula for the core for a new class of ideals, namely for the edge ideals  described above, Theorem~\ref{formulaext}.

The contents of this paper are as follows. We provide necessary
definitions and background material in Section~\ref{background}. In
Section~\ref{reductions} we discuss the format of minimal reductions
and restrictions on the coefficients  of their generators. In
Section~\ref{extension} we prove the main 
results of the paper,
namely that if $I$ is the edge ideal of any graph, then either $I$ is basic or $\core{I}\subset \m I$, Theorem~\ref{core in mI}, and if $I$ is the edge ideal of an even cycle with an arbitrary number of whiskers then 
$J:I=\m$ for every minimal reduction $J$ of $I$,
Theorem~\ref{colonm full whiskers}, and $\core{I}=\m I$,
Theorem~\ref{formulaext}. 
We give an example of a graph that is neither basic nor a whiskered even cycle for which  this formula for the core does not hold, see Example~\ref{counterexample}, and  the core is not a finite intersection of general minimal
  reductions. Furthermore,  Example~\ref{counterexample} 
  establishes that the condition that $I$ is weakly $(\ell
  -1)$-residually $S_2$ in \cite[Theorem~4.5]{CPU01} is
  necessary. 
  
  In general, the edge
  ideals of even cycles need not be weakly $(\ell-1)$-residually $S_2$. Therefore $\core{I}$ is not a priori a
  finite intersection of general minimal reductions in this
  case. Nevertheless,  in Section~\ref{finiteintersection} we show
  that the core of an even cycle is obtained via a finite intersection
  of 
  homogeneous binomial minimal reductions. It turns out these minimal
  binomial reductions also  establish the gradedcore. 
 We show that ${\rm graded}\core{I} = \core{I}$ for the edge ideals of even cycles, Remark~\ref{gradecore}.

\section{Background}\label{background}

Let $R$ be a Noetherian ring and $I$ an ideal. Suppose that
$I=(h_1, \ldots, h_q)$. The {\it Rees algebra} of $I$ is the subring
  $\mathcal{R}(I)=R[It] =R \oplus It
  \oplus I^2t^2 \oplus \ldots \subset R[t]$. 
 There is a canonical epimorphism $\phi: A=R[T_1, \ldots, T_q]  \longrightarrow \mathcal{R}(I)$
given by $T_i \mapsto h_it$. Let $L=\ker (\phi)$. Then $L= \bigoplus
\limits_{i=1}^{\infty}L_i $ is a graded ideal.
The ideal $I$ is said to be of {\it linear type} if
$L=L_1A$. It follows that  $J \subset I$ is a reduction of $I$ if and only if
$\mathcal{R}(I)$ is integral over $\mathcal{R}(J)$. Note that if $I$ is an ideal of linear
type then $I$ is basic.

Suppose $(R, \m,k)$ is a Noetherian local ring with infinite residue
field and $I$ is an ideal of $R$. The {\it special fiber ring} of $I$ is the
graded algebra \linebreak $\mathcal{F}(I)= \mathcal{R}(I) \otimes k= \bigoplus
\limits_{i\geq 0} I^{i}/\m I^{i}$. As above there is a canonical epimorphism
\linebreak $\psi: B=k[T_1, \ldots, T_q] \longrightarrow \mathcal{F}(I)$, 
whose kernel is a graded ideal referred to as the ideal of equations
of $\mathcal{F}(I)$. 

Northcott and Rees proved that when $R$ is a Noetherian local ring
then the minimal reductions 
correspond to Noether normalizations of $\mathcal{F}(I)$ \cite{NR}.  Furthermore, all minimal reductions have the same minimal
number of generators. This number is called the {\it analytic spread}
of $I$ and is defined by $\ell(I)=\dim \mathcal{F}(I)$. It then follows that $\mu(J)=\ell(I)$ for every minimal reduction $J$ of $I$
\cite{NR}. Throughout let $\ell=\ell(I)$ denote the analytic spread of
$I$.

Explicit descriptions of the Rees algebra, $\mathcal{R}(I)$,
and the special fiber ring  $\mathcal{F}(I)$  of an edge ideal $I$ were
obtained by Villarreal in \cite{Vill}. Let $G$ be a graph on a set of vertices $V=\{x_1, \ldots, x_n\}$. Define $I$ to be the ideal generated by
all elements of the form $x_ix_j$, where 
$\{x_i, x_j\}$ is an edge of $G$. Then $I=I(G)$ is the {\it edge
  ideal} associated to the graph $G$. In general, $I$ is an ideal of the
polynomial ring $k[x_1, \ldots, x_n]$ over a 
field $k$. As mentioned in Section~\ref{intro}, in order to discuss minimal
reductions of edge ideals of graphs, we will view $I$ as an ideal of the local
ring $R=k[x_1, \ldots, x_n]_{(x_1, \ldots,
  x_n)}$, where $k$ is an infinite field.
  
Villarreal characterized the edge ideals that are of linear type. More precisely, he showed that if $G$ is a
connected graph then the edge ideal of $G$ is of linear type if and
only if $G$ is a tree or has a unique cycle of odd length
\cite[Corollary~3.2]{Vill}. Since the edge ideals of odd cycles or trees are of linear
 type and 
hence have no proper reductions, these are precisely the graphs whose edge ideal is basic. Thus we will consider 
edge ideals of graphs with  {\it irreducible even closed walks}. Here a closed walk $x_1, e_1, x_2, e_2, x_3,\ldots ,e_d,x_1$ is considered to be reducible if there exists edges $e_i$ and $e_j$ in the walk such that $e_i$=$e_j$ and $i$ and $j$ have different parities. Such walks are considered reducible because they do not correspond to minimal relations of the defining ideal of the fiber cone \cite[Proposition~3.1]{Vill}. Note that a graph $G$ contains an irreducible even closed walk if and only if $G$ is not of linear type. Just as for a cycle, a closed walk is considered to be independent of its starting point for the purpose of uniqueness. This also allows an even closed walk to be represented by its edges with the vertices suppressed. Note that if $e_1, \ldots , e_d$ is an even closed walk, then $e_1, \ldots, e_d, e_1 \ldots , e_d$ is an even closed walk, which will be considered as a multiple of $e_1, \ldots , e_d$. A graph will be considered to have a unique irreducible even closed walk if all irreducible even closed walks are multiples of a fixed irreducible even closed walk.

Even cycles provide examples of irreducible even closed walks. For a more general example of an even closed walk, consider the graph whose edges are $e_1=x_1x_2, e_2=x_2x_3, e_3=x_1x_3, e_4=x_1x_4, e_5=x_4x_5, e_6=x_1x_5$. Then $ e_1,  e_2, e_3,  e_4, e_5,  e_6$ is an irreducible even closed walk without repeated edges that has a repeated vertex. For a nontrivial example of an irreducible even closed walk with repeated edges, consider the walk $e_1,e_2,e_3,e_4,e_5,e_6, \linebreak e_3,e_7$ in the graph whose edges are $e_1=x_1x_2, e_2=x_2x_3, e_3=x_3x_4, e_4=x_4x_5, e_5=x_5x_6, e_6=x_6x_4, e_7=x_3x_1$. Notice that if we label the edges of the walk $f_1, \ldots , f_8$, then $f_3=f_7$ and $3$, $7$ have the same parity.

\begin{remark}\label{basic properties}{\rm
Let $G$ be a graph with $s$ edges and a unique irreducible even closed walk  given by $e_{i_1}, e_{i_2}, \ldots, e_{i_d}$, and let $I=I(G)$ be the edge ideal of $G$. Then  $\mathcal{F}(I)\simeq k[T_1, T_2, \ldots,
T_s]/(T_{i_1}T_{i_3} \cdots T_{i_{d-1}}-T_{i_2}T_{i_4} \cdots T_{i_d})$, by
\cite[Proposition~3.1]{Vill}.  Therefore $\ell=s-1$ and  $I$ is an ideal of second analytic deviation one. 
}
\end{remark}

\section{The Structure of Minimal Reductions}\label{reductions}

We begin by proving a general result about the form of a minimal
reduction of an ideal $I$ of second analytic deviation one.  We state the following lemma for ease of reference.
 
\begin{lemma} $($\cite{NR}$)$ \label{NR red}
Let $(R,\m)$ be a Noetherian local ring. Let $I, K$ be ideals such
that $K \subset I$ and $\overline{K+\m I}=\overline{I}$, where
$\overline{I}$ denotes the integral closure of $I$. Then
$\overline{K}=\overline{I}$, i.e. $K$ is a reduction of $I$. 
\end{lemma}

\begin{lemma}\label{gen red format}
Let $R$ be a Noetherian local ring with infinite residue field. Assume
$I$ is an ideal with $\ell = \mu(I)-1$, and let $J$ be a
minimal reduction of $I$.  If $(h_1, \ldots, h_s)$ is a minimal
generating set of $I$, then $J$ has a generating set of
the form   
$(h_1+a_1h_{t}, h_2+a_2h_{t}, \ldots, h_{s}+a_{s}h_{t})$ 
for some $1 \leq t \leq s$, where $a_i \in R$ for all $i$ and $a_t=-1$.
\end{lemma}

\proof
Let $I=(h_1,\ldots, h_s)$ and let $J$ be a minimal reduction of
$I$.  If $s=1$ then the result is trivial. Suppose that $s \geq 2$.
Then $J=(f_1,\ldots, f_{s-1})$ for some $f_i \in I$. Let
$f_i=\sum \limits_{j=1}^{s}a_{ij}h_j$ and let $A=(a_{ij})$ be the
matrix of coefficients of $J$. Then $A$ is a $(s-1) \times s$
matrix. Let $\m$ denote the unique maximal ideal of $R$.

Suppose that $a_{ij} \in \m$ for all $i$ and $j$. Then $J \subset \m I
\subset I$. As $\ol{J}=\ol{I}$ then $\ol{0+\m I}=\ol{I}$. Hence by
Lemma~\ref{NR red} we have $0$ is a reduction of $I$, which is
impossible. Therefore $a_{ij} \not \in \m$ for some $a_{ij}$. After reordering the $h_i$ and
the $f_i$ we may assume, without loss of generality, that
$a_{11}=1$. Using row operations, which correspond to changing the
generating set of $J$, we can assume $A$
has the form  
$$\left(\begin{array}{ccccc} 1 & a_{12}  & \cdots  & a_{1,s-1} & a_{1,s} \\
0 & a_{22}  & \cdots  & a_{2,s-1} &a_{2,s} \\
\vdots & \vdots &  \ddots & \vdots &\vdots \\
0  & a_{s-1,2} & \cdots  & a_{s-1,s-1} & a_{s-1,s} \\
\end{array}\right).$$

Notice that $J$ is  minimally generated by $s-1$  elements (\cite{NR} or \cite[Proposition 8.3.7]{HS}). Hence the matrix $A$ has full rank and thus using an argument similar to the one above we may row reduce $A$ and assume that it is of the form
$$\left(\begin{array}{ccccc} 1 & 0  & \cdots & 0 & a_{1,s} \\
0 & 1  & \cdots  &0&a_{2,s} \\
\vdots & \vdots &  \ddots & \vdots &\vdots \\
0  & 0& \cdots  &1& a_{s-1,s} \\
\end{array}\right).$$
Then we may write $J$ as
$J=(h'_1+a_{1,s}h'_{s}, \ldots, h'_{i}+a_{i,s}h'_{s}, \ldots,
h'_{s-1}+a_{s-1,s}h'_{s})$, where $a_{i,s} \in R$ and $h_i'=h_{\sigma(i)}$
  for some permutation $\sigma$ of $\{1 \ldots, s\}$. The result follows
  by setting
  $t=\sigma (s)$, $a_t=-1$, and $a_{\sigma (i)}=a_{i,s}$ for all $1\leq i \leq s-1$.
\QED

The proof of Lemma~\ref{gen red format} can be extended for ideals with arbitrary second analytic deviation. 

\begin{corollary}\label{gen format}
Let $R$ be a Noetherian local ring with infinite residue field. Assume
$I$ is an ideal with $\ell = \mu(I)-n=s-n$, and let $J$ be a
minimal reduction of $I$.  If $(h_1, \ldots, h_s)$ is a minimal
generating set of $I$, then $J$ has a generating set of
the form   
$$(h_1+a_{1,1}h_{t_1}+\ldots +a_{1,n}h_{t_{n}}, \ldots, h_{s}+a_{s,1}h_{t_{1}}+\ldots +a_{s,n}h_{t_{n}})$$
for some $1 \leq t_1, \ldots , t_n \leq s$, where $a_{i,j} \in R$ for
all $i,j$ and $a_{t_{i},j}=-\delta_{ij}$  for all
$1 \leq i,j \leq n$.
 
\end{corollary}

Next we give an interpretation of Corollary~\ref{gen format} in the case of an edge ideal that contains a unique irreducible even closed walk.

\begin{corollary}\label{red format}
Let $I=(e_1, \ldots, e_s)$ be the edge ideal of a graph with $s$ edges containing a
unique irreducible even closed walk
and let $J$ be a 
minimal reduction of $I$. Then $J$ is of
the form   
$(e_1+a_1e_{t}, e_2+a_2e_{t}, \ldots, e_{s}+a_{s}e_{t})$
for some $1 \leq t \leq s$, where $a_i \in R$ for all $i$ and $a_t=-1$.
\end{corollary}

\proof
This follows immediately from Lemma~\ref{gen red format} and
Remark~\ref{basic properties}.
\QED

In addition to knowing the general form a reduction can take we also have control over the reduction number for the edge ideal of a graph with a unique irreducible even closed walk. 

Let $R$ be a Noetherian local ring, $I$ an ideal of $R$ and let $J$ be a minimal reduction of $I$.
The smallest $r$ for which the equality $I^{r+1}=JI^{r}$ holds is called the \textit{reduction
number of} $I$ \textit{with respect to} $J$ and is denoted by
$r_{J}(I)$. The reduction number $r_{J}(I)$ provides a measure of how closely
related $J$ is to $I$. The \textit{reduction
number} $r(I)$ of $I$ is the minimum of the reduction numbers
$r_{J}(I)$, where $J$ ranges over all minimal reductions of $I$. 
\begin{lemma}\label{red num}
Let $I$ be the edge ideal of a graph with $s$ edges containing a
unique irreducible even closed walk, which is of length $d$. Then  $r_J(I)= {\frac{d}{2}}-1$ for any minimal reduction $J$ of $I$.  In particular, $r_{J}(I)$ is independent of the minimal reduction $J$ of $I$.
\end{lemma}

\proof
By \cite[Proposition~3.1]{Vill} we know that the special fiber ring of $I$ is $\mathcal{F}(I)\simeq k[T_1, T_2, \ldots,
T_s]/(T_{i_1}T_{i_3} \cdots T_{i_{d-1}}-T_{i_2}T_{i_4} \cdots T_{i_d})$, where $e_{i_1}, \ldots ,e_{i_d}$ are the not necessarily distinct edges of the even walk.
Since the degree of the relation in the defining ideal of $\mathcal{F}(I)$ is ${\frac{d}{2}}$ then it follows that $r_J(I)= {\frac{d}{2}}-1$ by \cite[Proposition~5.1.3]{Vasc}. \QED

The next lemma and proposition allow us to use counting arguments to eliminate
potential reductions.

\begin{lemma}\label{K}
Let $I=(e_1, \ldots, e_s)$ be the edge ideal of a graph with $s$
edges, and let $J=(e_1+a_1e_{s},  \ldots, e_{s-1}+a_{s-1}e_{s})$. Fix $r
\geq 2$ and define 
$K^{r-1}$ to be the ideal generated by all elements of the form
$(e_i+a_ie_s)e_{i_{1}} \cdots e_{i_{r-1}}$ where $i \leq i_1 \leq i_2 \leq
\cdots \leq i_{r-1}$. Then $JI^{r-1}=K^{r-1}$. 
\end{lemma}

\proof 
For clarity, we first handle the case $r=2$. Clearly $K \subset
JI$. Since $JI$ can be generated by elements of the form $(e_q+a_qe_s)e_{i_1}$, we consider a generator  $(e_q+a_qe_s)e_{i_1}
\in JI$ for some $i_1<q<s$. 
Then
$$(e_q+a_qe_s)e_{i_1}=(e_{i_1}+a_{i_1}e_s)e_q - a_{i_1}(e_q+a_qe_s)e_s+
a_q(e_{i_1}+a_{i_1}e_s)e_s \in K.$$ 
Thus $JI=K$.

For the general case, consider a generator $(e_q+a_qe_s)M\in JI^{r-1}$,
where $M$ is a monomial generator of $I^{r-1}$. 
Write $M=e_{i_1}e_{i_2}\cdots e_{i_{r-1}}$ with $i_1 \leq i_2 \leq \cdots
\leq i_{r-1}$. Assume $i_1 < q$, and let $N=e_{i_2}e_{i_3}\cdots
e_{i_{r-1}}$.  
Then multiplying the equation above by $N$ yields
$$(e_q+a_qe_s)e_{i_1}N=(e_{i_1}+a_{i_1}e_s)Ne_q- a_{i_1}(e_q+a_qe_s)Ne_s +
a_q(e_{i_1}+a_{i_1}e_s)Ne_s.$$
Now by the choice of $i_1$, $(e_{i_1}+a_{i_1}e_s)Ne_q \in K^{r-1}$, as is
$(e_{i_1}+a_{i_1}e_s)Ne_s.$ Consider $(e_q+a_qe_s)Ne_s$. If $i_2 \geq q$
we are done. Otherwise, repeat the process for $Ne_s$. Since $M$ is a
product of $r-1$ edges, this process must terminate. Thus
$JI^{r-1} \subset K^{r-1}$. Since the other inclusion is clear,
$JI^{r-1}=K^{r-1}$ as claimed.
\qed

\begin{proposition}\label{numberofgenerators}
Let $I=(e_1, \ldots, e_s)$ be the edge ideal of a graph with $s$ edges
containing a 
unique irreducible even closed walk, which is of
length $d$. Let $J=(e_1+a_1e_{t},  \ldots, e_{s}+a_{s}e_{t})$  
for some $1 \leq t \leq s$, where
$a_i \in R$ and $a_t=-1$. Then   
$$\mu(I^r)=\left\{ \begin{array}{ll}
\binom{s+r-1}{r}, & r< {\frac{d}{2}} \\
\binom{s+r-1}{r}-1, & r=  {\frac{d}{2}}\end{array}\right.$$
and
$\mu(JI^{r-1}) \leq \binom{s+r-1}{r}-1$ for $r \geq 1$. 
\end{proposition}

\proof The number of products, allowing for 
repetition, of $r$ elements selected from a set containing $s$
elements is $\binom{s+r-1}{r}$, so $I^r$ can be generated by
$\binom{s+r-1}{r}$ monomials. From the structure of the fiber ring of
$I$, Remark~\ref{basic properties}, we
know  that there are no relations among the generators in degree less than
$\frac{d}{2}$, and there is precisely one relation in degree
$\frac{d}{2}$. 
Thus if $r<\frac{d}{2}$, there are no relations among the products
counted and the result 
follows. If $r=\frac{d}{2}$ and the edges of the 
irreducible even closed walk  are
$e_{i_1}, \ldots, e_{i_d}$, then $e_{i_1}e_{i_3}\cdots e_{i_{d-1}}=e_{i_2}e_{i_4}\cdots
e_{i_d}$ has been counted twice.  
Note that there are no other relations in degree
$\frac{d}{2}$ and  thus $\mu(I^r) = \binom{d+r-1}{r}-1$ for $r=\frac{d}{2}$.

Assume $J$ is an ideal of the given form. Select  any relabeling of the
edges of $G$ so that $t=s$. By Lemma~\ref{K}, in order to provide an upper
bound on the minimal number of 
generators of $JI^{r-1}$, it suffices to provide an upper bound on the
minimal number of generators of $K^{r-1}$. Note that for any $1 \leq i
<s$, there are $s-i+1$ generators of $I$ from
which $r-1$ are selected, with possible repetition, to form a monomial
$M$ for which $(e_i+a_ie_s)M$ is a generator of $K^{r-1}$.
There are $\binom{s-i+1+r-1-1}{r-1}$ possible
generators of $K^{r-1}$ of the form $(e_i+a_ie_s)M$ for each $1 \leq i
<s$. Now 
we have that $\sum \limits_{i=1}^s \binom{s+r-1-i}{r-1}=\binom{s+r-1}{r}$.
Thus there are 
$$\sum_{i=1}^{s-1} \binom{s+r-1-i}{r-1} = \binom{s+r-1}{r}-\binom{s+r-1-s}{r-1} =
\binom{s+r-1}{r} -1$$ 
elements in the generating set described above
for $K^{r-1}=JI^{r-1}$. This gives the desired upper bound on $\mu(JI^{r-1})$. 
\QED

Note that when $r < \frac{d}{2}$ the bound given above on the number of
generators of $JI^{r-1}$ is actually an equality. To see this, write $J=(f_1, \ldots , f_{s-1})$
and $I=(J,f_s)$ for some choice of $f_i$. Then among the generators
$f_{i_1} \cdots f_{i_r}$ of $I^r$, the only one that is not
automatically in $JI^{r-1}$ is $f_s^r$. Since
Proposition~\ref{numberofgenerators} shows that $I^r$ has
$\binom{s+r-1}{r}$ distinct generators for $r< {\frac{d}{2}}$, this
gives at least $\binom{s+r-1}{r}-1$ distinct generators of
$JI^{r-1}$. Thus if $r< {\frac{d}{2}}$ then $\mu(JI^{r-1}) =
\binom{s+r-1}{r}-1$. 

Using the information about the reduction numbers from Lemma~\ref{red
  num}  we show that the 
counting arguments used in Proposition~\ref{numberofgenerators} impose
restrictions on the coefficients of the generators of the reductions
in the case of edge ideals of graphs with a unique even cycle. Note
  that the proof below easily generalizes to graphs containing a
  unique even closed walk that does not contain repeated edges.    
Throughout the remainder of the paper, it will be convenient to reorder the edges of a cycle so that a particular edge is last. To that end, assume $e_1, \ldots , e_d$ form an even cycle, where $e_i=x_ix_{i+1}$ for $1 \leq i<d$ and
$e_d=x_1x_d$. We define a {\it cyclic
  reordering} of the vertices to be a relabeling $\sigma$ of the
vertices such that $\sigma(x_i) = x_{i+j}$ for some fixed $j$, where
subscripts are taken modulo $d$ and $\ol{0}=d$. Such a reordering preserves
adjacencies and the cycle structure, but allows any particular edge of
the cycle to be considered last, namely as  $e_d$.

\begin{corollary}\label{productofbs} 
Let $I=(e_1, \ldots , e_s)$ be the edge ideal of a graph with $s$
edges containing a 
unique even cycle, $e_1, \ldots ,e_d$. Define $J=(e_1+a_1e_{t},  \ldots, e_{s}+a_{s}e_{t})$  
for some $1 \leq t \leq s$, where
$a_i \in R$ and $a_t=-1$.  If $\prod \limits_{i=1}^{{\frac{d}{2}}}
a_{2i-1} = \prod 
  \limits_{j=1}^{{\frac{d}{2}}} a_{2j}$, 
then $J$ is not a reduction of $I$. 
\end{corollary}

\proof
If $J$ is a reduction of $I$, then $J$ must be minimal since it has
$\ell$ generators. By Lemma~\ref{red num}, we know that $J$ is
a minimal reduction of $I$ if 
and only if $JI^{r-1} =I^r$, where $r=\frac{d}{2}$. 

There are two cases to consider. 
If $t \leq d$, then after a cyclic
reordering of the cycle we may assume $t=d$ and $a_d=-1$. Otherwise, $t>d$.
Assume $\prod \limits_{i=1}^{{\frac{d}{2}}} a_{2i-1} = \prod
  \limits_{j=1}^{{\frac{d}{2}}} a_{2j}$.  Using this equality and the
  relation among the edges of the cycle, it is 
easy to check that for $t \geq d$
$$\begin{array}{lll}&&(e_1+a_1e_t)e_3e_5\cdots e_{d-1}=\\
&=& \sum  \limits_{i=1}^{r} (-1)^{i-1} a_2a_4 \cdots
  a_{2i-2}(e_{2i}+a_{2i}e_t)e_t^{i-1}e_d e_{d-2} \cdots e_{2i+2}\\ 
&+&  \sum  \limits_{i=1}^{r-1} (-1)^{i-1} a_1a_3 \cdots
 a_{2i-1}(e_{2i+1}+a_{2i+1}e_t)e_t^{i}e_{d-1} \cdots e_{2i+3}, 
\end{array}$$
where empty products are defined to be one. Note that this is a relation among the generators of $K^{r-1}$ that were counted in Proposition~\ref{numberofgenerators}.  Therefore by Lemma~\ref{K}, $\mu(JI^{r-1})=\mu(K^{r-1}) \leq  \binom{d+r-1}{r}-1-1  < \mu(I^{r})$. Thus $J$ is
not a reduction of $I$. \QED

We conclude this section by providing
concrete examples of reductions for the edge ideals of graphs containing a
unique irreducible even closed walk. Note that these examples will  provide the building blocks for computing the
core as a finite 
intersection in Section~\ref{finiteintersection}.

\begin{example}\label{basicreduction}
Let $I$ be the edge ideal of a graph of an even cycle. Let $R$ be the corresponding localized polynomial ring and let $k$ be the residue field of $R$. We further assume that 
the characteristic of $k$ is not $2$. Let $J=(e_1+a_1e_{t},  \ldots, e_{d}+a_{d}e_{t})$
for some $1 \leq t \leq d$, where
$a_i=1$ for all $i \neq t$ and $a_t=-1$. Then $J$  is a minimal reduction of $I$.
\end{example}

\proof
If $J$ is a reduction of $I$, then $J$ is a minimal reduction since $J \subset I$ and $\mu(J)=\ell$. After a cyclic reordering we may assume $t=d$ and $a_d=-1$. Let
$r={\frac{d}{2}}$. Clearly $JI^{r-1} \subset I^r$. To see the other
  inclusion, we first prove $e_d^r \in JI^{r-1}$. Notice that
$e_d^r + (-1)^{r-1} \prod \limits_{i=1}^{r} e_{2i-1} \in JI^{r-1}$
since 
$$\begin{array}{lll}
e_d^r+ (-1)^{r-1} \prod \limits_{i=1}^{r} e_{2i-1}&=& 
\sum \limits_{i=1}^{r}(-1)^{i-1} (e_{2i-1}+e_d)e_1 \cdots e_{2i-3}e_d^{r-i},
  \end{array}$$
where empty products are defined to be one.
Similarly, 
$e_d^r + (-1)^{r-2} \prod \limits_{j=1}^{r} e_{2j} \in JI^{r-1}$
since $ e_d^r+ (-1)^{r-2} \prod \limits_{j=1}^{r} e_{2j}=  \sum
\limits_{i=1}^{r-1}(-1)^{i-1} (e_{2i}+e_d)e_2 \cdots
e_{2i-2}e_d^{r-i}$.
Combining these relations with the relation on the edges $\prod \limits_{i=1}^r
e_{2i-1}=\prod \limits_{j=1}^r e_{2j}$ 
gives $2e_d^r \in JI^{r-1}$. Thus $e_d^r \in JI^{r-1}$ as desired.

Now let $M\in I^r$ be a monomial generator. If $M=e_d^r$ we are done
by the argument above. If not, write $M=e_{i_1}e_{i_2}\cdots e_{i_r}$
for some choice of $r$ edges, ordered so that $i_1 \leq i_2 \leq
\cdots \leq i_r$. Define $M_1=e_{i_2}e_{i_3}\cdots
e_{i_r}$ and consider $(e_{i_{1}}+e_d)M_1=M+e_dM_1.$ If $M_1=e_d^{r-1}$,
then since $e_dM_1$ and $(e_{i_{1}}+e_d)M_1$ are both in $JI^{r-1}$, we see
that $M \in JI^{r-1}$ as well. If $M_1 \not= e_d^{r-1}$, then
define $M_2=e_{i_3}e_{i_4}\cdots e_{i_r}$. Notice
that if $M_2=e_d^{r-2}$, then by the equation
$(e_{i_2}+e_d)M_2=M_1+e_dM_2$ one sees that $M_1\in JI^{r-2}$ as above,
which then implies $M\in JI^{r-1}$. If $M_2 \not= e_d^{r-2}$ we
repeat the process. The process is clearly finite, and since at each
stage of the algorithm, $M_i$ is replaced by $e_dM_{i+1}$, the
algorithm will terminate. Thus for some (not
necessarily distinct) edges  $e_{i_j}$,
$M+(-1)^{q-1}e_d^r=(e_{i_1}+e_d)M_1-(e_{i_2}+e_d)e_dM_2+\ldots+
(-1)^{q-1}(e_{i_q}+e_d)e_d^{q-1}M_q$, where $q \leq {\frac{d}{2}}$ and
$M_q=e_d^{r-q}$. Thus $M \in JI^{r-1}$.
\QED

Example~\ref{basicreduction} generalizes to even closed walks without repeated edges.
We remark that when ${\rm char}\; k=2$ then it follows immediately from
Corollary~\ref{productofbs} that the ideal $J$ in
Example~\ref{basicreduction} is not a minimal reduction of $I$. In
order to avoid  
characteristic dependent arguments,  we provide 
two additional examples of minimal reductions that are free of
characteristic assumptions and which hold for edge ideals of graphs containing a (not necessarily unique) irreducible even
  closed walk.

\begin{example}\label{oddreduction}
Let $I$ be the edge ideal of a graph containing an irreducible 
even closed walk $e_1, \ldots, e_d$. Write $I=(e_1, \ldots ,e_d,e_{d+1}, \ldots, e_s)$, where  $e_{d+1}, \linebreak \ldots ,e_s$ are the distinct edges of $G$ not contained in the walk. Define $\delta'_{i,j}$ to be $-1$ if $e_i=e_j$ and $1$ otherwise.
Then $J=(e_1, e_2+\delta'_{2,d}e_d,e_3,
\ldots ,e_{d-2}+\delta'_{d-2,d}e_d, e_{d-1}, e_{d+1}+e_d, \ldots , e_s+e_d)$ is a reduction of $I$.  Furthermore, if $I$ contains a unique irreducible even closed walk, then $J$ is a minimal reduction of $I$.
\end{example}

\proof
Note that the first $d$ generators of $I$ are not necessarily unique, but that any repeated edges will have the same parity. Also, any repeated edge other than $e_d$ listed in the generating set of $I$ corresponds to a repeated generator of $J$. Hence $\mu(J)=\mu(I)-1$. Let $r={\frac{d}{2}}$. Clearly $JI^{r-1} \subset I^r$. For the other
inclusion,  let $M$ be a monomial generator of $I^{r}$. 
Write $M=e_{i_1}e_{i_2}\cdots e_{i_r}$
for some choice of $r$ edges, where if a repeated edge divides $M$, the largest possible subscript for the edge is used.
If $i_j$ is odd and less than $d$ for some $j$, then $M=e_{i_j}N$,
where $e_{i_j} \in J$ 
and $N\in I^{r-1}$. Thus $M\in JI^{r-1}$. So suppose $i_j$ is not odd
for all $i_j<d$. Define $s_i$ to be the number of times that $e_{j}=e_d$ for $j < i$. As in
Example~\ref{basicreduction} we have 
$$e_d^r + (-1)^{r-2-s_d} \prod
\limits_{j=1}^{r} e_{2j} =  \sum
\limits_{i=1}^{r-1}(-1)^{i-1-s_{2i}} (e_{2i}+\delta'_{2i,d}e_d)e_2 \cdots
e_{2i-2}e_d^{r-i} \in JI^{r-1}.$$ 
By the relation $\prod \limits_{i=1}^r
e_{2i-1}=\prod \limits_{j=1}^r e_{2j}$ and the fact that $\prod
\limits_{i=1}^r e_{2i-1} \in JI^{r-1}$ we have that $e_{d}^{r} \in JI^{r-1}$.
The remainder of the argument follows as in
Example~\ref{basicreduction} by noting that each $e_{i_j}$ in the
expression for 
$M$ now has $i_j$ even or $i_j\geq d$ and thus $(e_{i_j}+e_d)M_j \in JI^{r-1}$ for
each $j$. Finally, when $I$ contains a unique irreducible even closed walk  then $\ell=\mu(I)-1$.  Hence $J$ is a minimal reduction of $I$. \qed

\begin{example}\label{evenreduction}
Let $I$ be the edge ideal of a graph containing an irreducible 
even closed walk $e_1, \ldots, e_d$. Write $I=(e_1, \ldots ,e_d,e_{d+1}, \ldots, e_s)$, where  $e_{d+1},\linebreak \ldots ,e_s$ are the distinct edges of $G$ not contained in the walk. Define $\delta_{i,j}$ to be $0$ if $e_i=e_j$ and $1$ otherwise.
Then $J=(e_1+e_d,\delta_{2,d}e_2,e_3+e_{d}, \ldots,
\delta_{d-2,d}e_{d-2} , \linebreak e_{d-1}+e_d, e_{d+1}+e_d,  \ldots , e_s+e_d)$ is a reduction of $I$. Furthermore, if $I$ contains a unique irreducible even closed walk, then $J$ is a minimal reduction of $I$.\end{example}

\proof
The proof is similar to the proof of Example~\ref{oddreduction}.
\qed

\section{Cores of Edge Ideals of Whiskered Cycles}\label{extension}

Recall that if $I$ is the edge ideal of a connected graph,
then $I$ is of linear type if and only if $I$ is the edge ideal of a
tree or of a graph 
containing a unique cycle of odd length  by
\cite[Corollary~3.2]{Vill}, and thus
$\core{I}=I$.  This implies that $I$ is not of linear type if and only
if the graph associated to $I$ has an irreducible even closed walk. 
In this section, we show that if
$I$ is the edge ideal of any graph that is not basic, then we have
$\core{I} \subset \f{m}I$. We also establish a class of graphs for
which this inclusion is an equality. Note that the core of a monomial ideal
is also a monomial ideal by \cite[Remark~5.1]{CPU01}.

\begin{theorem}\label{core in mI}
Let $I$ be the edge ideal of a connected graph containing an
irreducible even closed walk. Then
$\core{I} \subset \f{m}I$.
\end{theorem}

\proof 
Write $I=(e_1, \ldots ,e_s)$, where $e_1, \ldots ,e_d$ form an irreducible even closed walk. Let $e_i$ be a generator of $I$.  If $i$ is odd then
$$J_1= (e_1+e_d,\delta_{2,d}e_2, \ldots,
\delta_{d-2,d}e_{d-2} ,e_{d-1}+e_d, e_{d+1}+e_d,  \ldots , e_s+e_d)$$ is a reduction of $I$ by
Example~\ref{evenreduction} and $e_i \notin J_1$. Similarly, if $i$ is
even then $J_2=(e_1, e_2+\delta'_{2,d}e_d,
\ldots ,e_{d-2}+\delta'_{d-2,d}e_d, e_{d-1}, e_{d+1}+e_d, \ldots , e_s+e_d)$ is a
reduction of $I$ by Example~\ref{oddreduction} and $e_i \notin J_2$. Therefore $e_i \not\in
\core{I}$.

Let $g$ be a minimal monomial generator of $\core{I}$. Since $g \in I$  then $g=fe_i$ for some
$e_i$ and $f\in R$ a monomial. Since $e_i \not \in \core{I}$ then $f
\in \f{m}$. Therefore $g \in \f{m}I$ and thus $\core{I} \subset
\f{m}I$. 
\QED

We state the following result without a proof, as its proof is elementary.

\begin{lemma}\label{det B}
Let $R$ be a commutative ring with identity, let $d \geq 4$ be an even
integer, and let  $b_1, \ldots,
b_{d} \in R$. Let $B$ be a $d\times d$ matrix of the following form: 
$$
B= \left(\begin{array}{ccccccc}
0 & b_d & 0 & 0 & \ldots & 0 & -b_1\\
-b_2 & 0 & b_1 & 0 &0 & \ldots & 0\\
0 & -b_3 & 0 & b_2 & 0 & \ldots & 0\\
\vdots & \vdots & \vdots & \vdots & \vdots & \vdots & \vdots \\
b_{d-1} & 0 & 0 & \ldots & 0 & -b_d & 0
\end{array}\right).$$
 Then $\det B =\left( \prod \limits_{i=1}^{\frac{d}{2}} b_{2i-1} -  \prod
\limits_{j=1}^{\frac{d}{2}} b_{2j} \right)^2$.
\end{lemma}

For the rest of the article we will assume that $I$ is the edge ideal of a graph $G$ with a unique even cycle and will 
order the edges so that $e_1, \ldots, e_s$ are the edges of $G$ and
$e_1, \ldots , e_d$ are the edges of the even cycle. In general, if $G$ is a connected graph on $n$ vertices with $s$
edges, then $s \geq n$ with equality if and only if $G$ has a unique 
cycle. Thus for the remainder of the article, the number of edges will be the same as the number of vertices of the graph. For the next
theorem, we need to further restrict the class of graphs
considered. 

\begin{assumptions}\label{assm even whiskers}{\rm 
Let $G$ be a connected graph on the vertices $x_1, \ldots , x_s$ containing a
unique cycle, which is of even length $d\geq 4$, given 
by $e_i=x_ix_{i+1}$ for $1 \leq i<d$ and
$e_d=x_1x_d$. Assume
further that $x_j$ is a leaf for all $j >d$. Thus for each $ j>d$
there exists a unique vertex $x_{i_{j}}$ with $1 \leq i_j \leq d$  such that $e_j=x_{i_j}x_j$
is an edge of $G$.
Notice that  it is not required that the $i_j$
be distinct for different $j$. It is possible for a single vertex of
the cycle to have multiple leaves as neighbors.
Let $I=(e_1,\ldots, e_{s})$ be the edge ideal of $G$ in the localized
polynomial ring  $R=k[x_1, \ldots, x_s]_{(x_1, \ldots, x_s)}$ over an
infinite field $k$. Then $\mu(I)=s$,  and $\ell=s-1$ by
\cite[Proposition~3.1]{Vill}. We remark that Corollary~\ref {productofbs}  holds for this class of ideals.

}  
\end{assumptions}

The following theorem shows that for the class of  edge ideals $I$ with a unique even cycle and an arbitrary number of whiskers, the ideal $J:I$ is independent of the minimal reduction
$J$ of $I$.

\begin{theorem} \label{colonm full whiskers}
Let $R$ and $I$ be as in \ref{assm even whiskers}  and let $J$ be a minimal reduction of $I$. Then $J:I=\m$.

\end{theorem}

\proof 
Let $J$ be a minimal reduction of $I$. Then  $J$ is of the form $(e_1+a_1e_{t}, \ldots, e_{s}+a_{s}e_{t})$  for some $1 \leq t \leq s$, where $a_i \in R$ for all $i$ and $a_t=-1$, by Corollary~\ref{red format}.
Let $f_i=e_i+b_ie_t$, where $b_i=a_i$ if $a_i \not \in
\m$ and $b_i=0$ if $a_i \in \m$. Consider $J'=(f_1, \ldots, f_{s})$, where $f_t=0$ since $b_t=-1$.
Notice that $J \subset J'+\m I \subset I$. Hence $J'$ is a
reduction of $I$ by Lemma~\ref{NR red}.  

Consider a presentation matrix $\phi$ of $I$, where $R^{q}
\stackrel{\phi}{\longrightarrow} R^{s} \longrightarrow I
\longrightarrow 0$. Let $\psi$ be the submatrix of $\phi$ consisting of the
linear relations on the generators of $I$. Then $\psi$ is an $s \times (2s-d)$ matrix of the form $\psi =\left(\begin{array}{ccc}
\psi_1 & \psi_2 &\psi_3\\
\end{array}\right)$, where $\psi_1, \psi_2, \psi_3$ are matrices
defined below. For the remainder of the proof we let ${\ol{ i}}=i$ modulo $d$, with the convention that $\ol{0}=d$.

Let $\psi_1$ be an $s \times d$ matrix such that for each $ 1 \leq i
\leq d$ the $i$-th column is $(0, \ldots, 0,-x_{\ol{i+1}}, x_{\ol{ i-1}}, 0,
\ldots, 0)^{T}$, where $-x_{\ol{ i+1}}$ is the 
${\ol{ (i-1)}}$ entry and $x_{\ol{ i-1}}$ is the $i$-th entry.

Let $\psi_2$ be an $s \times (s-d)$ matrix such that  for each $
d+1\leq j \leq s$ the $(j-d)$-th column is $(0,  \ldots, 0, x_j,  0,
\ldots, 0, -x_{\ol{ i_{j}-1}}, 0,\ldots, 0)^{T}$, where $x_j$ is the
${\ol{ (i_j-1)}}$ entry and $ -x_{\ol{ i_{j}-1}}$ is the $j$-th entry.

Let $\psi_3$ be an $s \times (s-d)$ matrix such that for each $
d+1\leq j \leq s$ the $(j-d)$-th column is $(0,  \ldots, 0, x_j,  0,
\ldots, 0, -x_{\ol{ i_{j}+1}}, 0,\ldots, 0)^{T}$, where $x_j$ is the
$i_j$ entry and $ -x_{\ol{i_{j}+1}}$ is the $j$-th entry. 

We remark that if $s=d$, then the matrices $\psi_2$ and $\psi_3$ are zero and the matrix $\psi$ is a $d \times d$ matrix. 
Notice that performing a
series of elementary row operations on $\phi$ corresponds to altering
the generating set of $I$. We choose elementary row operations so that the
generating set of $I$ becomes $I=(J', e_t)$. Let
$\phi'$ be the corresponding presentation matrix of $I$ and $\psi'$ the submatrix
consisting of the columns containing the linear relations.
By the choice of the generating set, the $t$-th row of $\phi'$ forms
a (not necessarily minimal) presentation matrix $\tilde{\phi}$ of
$I/J'$.  Let $\widetilde{\psi}$ denote the $t$-th row of $\psi'$. We will show that $I_1(\widetilde{\psi})=I_1(\tilde{\phi})=\m$. Notice that

$ \begin{array}{lll} I_{1}(\widetilde{\psi})&=&(\{b_{\ol{i-1}}x_{\ol{i+1}}-b_{i}x_{\ol{i-1}} \mid \mbox{ for }1\leq i \leq d\}, \\
&&\{b_jx_{\ol{ i_{j}-1}}-b_{\ol{ i_{j}-1}}x_j \mid \mbox{ for } d+1 \leq j \leq s\}, \\
 &&\{b_jx_{\ol{ i_{j}+1}}-b_{i_{j}}x_j \mid \mbox{ for } d+1 \leq j \leq s\} ). \end{array}$

Then $\widetilde{\psi}^T = B \cdot (\underline{x})^T$, where
$(\underline{x})^{T}=(  x_1, \ldots, x_s)^{T} $ and
$B=\left(\begin{array}{c} B_1  \\ 
C \end{array} \right)$, where $B_1=\left(\begin{array}{cc} B_0 &
  \underline{0} \end{array} \right)$, 

$B_0=\left(\begin{array}{ccccccc}
0 & b_d & 0 & 0 & \ldots & 0 & -b_1\\
-b_2 & 0 & b_1 & 0 &0 & \ldots & 0\\
0 & -b_3 & 0 & b_2 & 0 & \ldots & 0\\
\vdots & \vdots & \vdots & \vdots & \vdots & \vdots & \vdots  \\
b_{d-1} & 0 & 0 & \ldots & 0 & -b_d & 0 
\end{array}\right),$

\noindent  $\underline{0}$ is the $d \times (s-d)$ zero matrix, and  $C$ is a 
$2(s-d) \times s$ matrix. We construct $C$ as follows. For each $d+1 \leq j \leq s$ there are two rows of $C$:

$C_{j}^{1}=(0, \ldots, 0, b_j, 0, \ldots,0, -b_{\ol{ i_{j}-1}}, 0, \ldots,
0)$, where $b_j$ is the ${\ol {(i_j-1)}}$-th entry and $ -b_{\ol{ i_{j}-1}}$ is the $j$-th entry and 

$C_{j}^{2}=(0, \ldots, 0, b_j, 0, \ldots,0, -b_{i_{j}}, 0, \ldots, 0)$,
where $b_j$ is the ${\ol{ (i_j+1)}}$-th entry and $ -b_{i_{j}}$ is the $j$-th entry.

Notice that when $s=d$ then $B=B_0$ and $\det B=\det B_0 \neq 0$, by Lemma~\ref{det B} and Corollary~\ref{productofbs}. In general, we will construct an $s \times s$ submatrix of $B$ with a nonzero determinant and thus after row reducing $B$ we will have
$I_{1}(\widetilde{\psi})= \m$.

We remark that by construction of the submatrix $C$, for each $d+1
\leq j \leq s$ the rows $C_j^1$ and $C_j^2$ have nonzero entries in the $j$-th column,
one of those entries is $-b_{\ol{ i_{j}-1}}$ and the other is $-b_{i_j}$. Notice that one of
$\ol{{i_{j}-1}}$ and ${i_j}$ will be even and one will be odd.

First consider the submatrix $B_1'=\left(\begin{array}{c} B_1  \\
C_{1}   \end{array} \right)$, where $C_{1}$ is the submatrix of $C$ constructed by 
selecting all the rows  of $C$ such that for each ${d+1}~\leq j \leq
s$ the entry in the $j$-th column is $-b_r$ for some  $r$
even. Notice that $B_1'$ is a block matrix and after exchanging rows
of $C_1$ we have a diagonal matrix of size $(s-d) \times (s-d)$ in the lower right corner. Thus  after these row operations $B_1'$ is equivalent to $\left(\begin{array}{cc} B_0 & \underline{0} \\
C_{1}' & D_1  \end{array} \right)$, where $D_1$ is diagonal with
diagonal entries of the form $-b_r$ with $ 2 \leq r \leq d$
even.  

Therefore, $\det B_1'=\pm \det B_0 \det D_1$.
Since $D_1$ is diagonal, $\det D_1$ is the product of its diagonal
entries. Notice that each diagonal entry of $D_1$ is by definition of
the form $b_r$ for some even $2 \leq r \leq d$, but not all even
$r$ need occur, and some could occur multiple times. 

We now consider another $s \times s$ submatrix of $B$, namely
$B_2'=\left(\begin{array}{c} B_1   \\ 
C_{2}   \end{array} \right)$, where $C_2$ is the submatrix of $C$ constructed by 
selecting all the rows  of $C$ such that for each  $d+1 \leq j \leq s$
the entry in the $j$-th column is $-b_q$ such that $q$ is odd. Notice
that $B_2'$ is a block matrix and after exchanging rows of $C_2$ we
have a diagonal matrix of size $(s-d) \times (s-d)$ in the lower
right corner. Thus $B_2'$ is equivalent to $\left(\begin{array}{cc} B_0 &
  \underline{0} \\ 
C_{2}' & D_2  \end{array} \right)$, where $D_2$ is diagonal with
diagonal entries of the form $-b_q$ with $ 1 \leq q \leq d$
odd. Notice that the diagonal entries of $D_2$ are not necessarily
distinct. 
As before $\det B_2'= \pm \det B_0 \det D_2$  and  $\det D_2$ is a
product of its diagonal entries, each of which has an odd subscript.

We observe that  $\det B_1'$ and $\det B_2'$ are not simultaneously
zero. By Corollary~\ref{productofbs} and Lemma~\ref{det B} we have $\det B_0=\left( \prod \limits_{i=1}^{\frac{d}{2}} b_{2i-1} -  \prod\limits_{j=1}^{\frac{d}{2}} b_{2j} \right)^2 \neq
0$. It follows that since each $b_i \in k$ it is not
possible to have $b_q=0$ for some odd $q$ and $b_r=0$ for some even
$r$ simultaneously. Thus $\det D_1$ and $\det D_2$ cannot be
simultaneously zero.

Therefore $I_{1}(\widetilde{\psi})= \m$. 
Notice that
we have 
$R^{q} \stackrel{\tilde{\phi}}{\longrightarrow} R \longrightarrow
I/J' \longrightarrow 0$ and $I_{1}(\widetilde{\psi}) \subset
I_1(\tilde{\phi}) \subset {\rm ann} (I/J')=J':I$. Furthermore, since
$J' $ is a minimal reduction of $I$ then $J':I \neq R$. Hence
$I_1(\tilde{\phi})=\m=J':I$.  

Recall that $J \subset J' +\m I \subset I$. Since $J':I= \m$ then $\m I \subset J'$ and thus $J \subset J'
\subset I$. Since $J$ and $J'$ are both minimal reductions of $I$ and $J \subset J'$ then $J=J'$ and thus $J:I=\m$ as well.\QED

A careful examination of the above proof shows that it yields even more
information about the form a minimal reduction can take. In particular, the
coefficients $a_i$ of Corollary~\ref{red format} can be taken to be units.

\begin{corollary}\label{unit coeff full whiskers}
Let $R$ and $I$ be as in ~\ref{assm even whiskers}, and let $J$ be a
minimal reduction of $I$.  Then $J$ is  of the form   
$(e_1+b_1e_{t}, \ldots, e_{t}+b_te_{t}, \ldots, e_{s}+b_{s}e_{t})$ for some $t$,
where $b_t=-1$ and for $1\leq i \leq s$, either $b_i \not\in \m$ or
$b_i=0$.

\end{corollary}
\proof By Corollary~\ref{red format} there exist $a_i \in R$ such that $J=(e_1+a_1e_{t},
 \ldots, e_{i}+a_ie_{t}, \ldots, e_{s}+a_{s}e_{t})$, where $a_t=-1$. 
 Let $b_i=a_i$ if $a_i \not\in \m$ and $b_i=0$ if $a_i \in \m$. Then by the proof
 of Theorem~\ref{colonm full whiskers} we have that $J=J'=(e_1+b_1e_{t}, \ldots,
 e_{i}+b_ie_{d}, \ldots, e_{s}+b_{s}e_{t})$. \QED

We are now ready to prove the second  main theorem of this section. 

\begin{theorem} \label{formulaext}
Let $R$ and $I$ be as in \ref{assm even whiskers}. Then
$\core{I}=(J:I)I=\m I$
for any minimal reduction $J$ of $I$.
\end{theorem}

\proof
By Theorem~\ref{colonm full whiskers} we have 
$J:I=\m$ for every minimal reduction $J$ of $I$. Hence for any minimal
reductions $J$ and $J'$ of $I$ we have $J:I=J':I$. In particular, $(J:I)I \subset J'$ and thus $\m
I=(J:I)I \subset \core{I}$. By Theorem~\ref{core in mI}
we have the other inclusion and thus $\core{I}=(J:I)I=\m I$. \hfill \QED

\begin{remark}{\rm
Let $R$ be a Gorenstein local ring 
  with infinite residue field and $I$ an ideal that satisfies $\depth R/I^{j} \geq \dim R/I -j+1$ for
all $1 \leq j \leq \ell-g+1$, where $g={\rm ht}\; I>0$. We further assume that $I$ satisfies  $G_{\ell}$. This condition is rather mild; it requires that $\mu(I_{\f{p}}) \leq \dim R_{\f{p}}$ for every prime $\f{p}$ containing $I$ with $\dim R_{\f{p}} \leq \ell-1$.  Under these assumptions $r(I) \leq
\ell-g+1$ is equivalent to $\core{I}=(J:I)J=(J:I)I$
for every minimal reduction $J$ of $I$  as was shown in \cite[Theorem~2.6, 
  Corollary~3.7]{CPU01}. Therefore the formula for the core we obtain in Theorem~\ref{formulaext} is not surprising. We remark that edge ideals of even cycles do  satisfy $G_{\ell}$ but the depth condition  above does not hold for the edge
ideals of even cycles of length $d \geq 6$ and thus our result does not follow from \cite[Theorem~2.6]{CPU01}. Nonetheless the reduction number
  for these ideals 
is $r(I)=\frac{d}{2}=\ell-g <\ell-g+1$ as shown in Lemma~\ref{red num}. 
}
\end{remark}

Before we can proceed  we need to recall some
definitions. Let $R$ be a Noetherian ring and $I$ an ideal of ${\rm
  ht} \; I=g >0$. For each $ i \geq g$ a {\it geometric $i$-residual
  intersection} of $I$ is an ideal $K$ such that there exists an
$i$-generated ideal $\mathfrak{a} \subset I$ with $K=\mathfrak{a}:I$, ${\rm
  ht}\; K \geq i$, and ${\rm ht} (I+K) \geq i+1$. 
Furthermore, $I$ is {\it weakly $n$-residually $S_2$} if $R/K$ satisfies Serre's condition $S_2$
for every geometric $i$-residual intersection $K$ of $I$ and for all $g \leq i \leq n$.

The following example shows that the formula
for the core given in Theorem~\ref{formulaext} does not hold in general if $I$ is the edge ideal of a
graph with a unique cycle that is even. 

\begin{example}\label{counterexample}{\rm 
Let $G$ be a graph on the vertices $x_1,  \ldots, x_6$ with
edges  \linebreak $e_1=x_1x_2, e_2=x_2x_3, e_3=x_3x_4, e_4=x_1x_4,e_5=x_4x_5,e_6=x_5x_6$.
Let $I$ be the
edge ideal of $G$ in $R={\mathbb Q}[x_1, \ldots, x_6]_{(x_1, \ldots, x_6)}$
and let
$\m=(x_1, \ldots, x_6)$ denote the maximal ideal of $R$. Then $  \m I
\not\subset \core{I} $. Furthermore, $I$ is not weakly \linebreak $(\ell -1)$-residually $S_2$ and $\core{I}$ is not a finite intersection of general minimal reductions of $I$.

 }

\end{example}

\proof
Notice that the graph $G$ is a square with two additional edges. By Remark~\ref{basic properties} we know that $\ell=5$. Also $g={\rm ht}\;I=3$.
Let  $H=(e_1+e_2,e_3+e_2,e_4+e_2,e_5,e_6+e_2)$. 
It is straightforward to verify that $I^2=HI$ and thus $H$ is a
minimal reduction of $I$. Using Macaulay
2~\cite{M2} we see that  $H:I=(x_1, \ldots, x_5)$.  Therefore, if $\m
I \subset \core{I}$ then $\m I \subset H$ and thus $\m \subset H:I$, a
contradiction.  Hence $\m I \not\subset \core{I}$

We will now show that $\core{I}$ is not a finite intersection of general minimal reductions of $I$.
We follow the outline of the proof of Theorem~\ref{colonm full whiskers}. Let $\phi$ be a presentation matrix of $I$. Then
the matrix $\psi$ of the linear relations on the generators of $I$ is given by 
$$\psi=\left(\begin{array}{ccccccc} x_4 & -x_3 &0 & 0 & 0 &0 &0\\
0 & x_1 &-x_4 & 0 &0 &0&0 \\
0 & 0& x_2  & -x_1   &0&x_5 &0 \\
-x_2 & 0& 0  & x_3   &-x_5&0 &0 \\
0  & 0&0  &0  &x_1&-x_3& x_{6} \\
0  & 0&0  &0  &0&0&-x_{4} \\
\end{array}\right).$$

Let $J$ be a minimal reduction of $I$. Then by Corollary~\ref{red format} we obtain that  $J=(e_1+a_1e_t, \ldots, e_6+a_6e_t)$,
where $ 1 \leq t \leq 6$, $a_t=-1$, and $a_j \in R$ for all $1 \leq j
\leq 6$. Let  $f_j=e_j+b_je_t$, where $b_j=a_j$ if $a_j \not\in \m$
and $b_j=0$ if $a_j \in \m$ for $1\leq j \leq 6$. Let $J'=(f_1, \ldots
f_6)$. Notice that $f_t=0$ since $b_t=-1$, and $J \subset J'+\m I \subset
I$. Therefore $J'$ is also a reduction of $I$ by Lemma~\ref{NR red}.
Then $I=(J, e_t)$. We choose
elementary row operations so that $\phi'$ is the new presentation matrix
of $I$ that reflects the
 generating set $(J,e_t)$ of $I$ and $\psi'$ is the corresponding matrix of linear relations. Notice that by the choice of the generating
set for $I$, the $t$-th row of $I$ forms a (not necessarily minimal)
presentation 
matrix $\widetilde{\psi}$ of $I/J'$. Then
$I_1(\widetilde{\psi})=(b_4x_2-b_1x_4,b_1x_3-b_2x_1, b_2x_4-b_3x_2,
b_3x_1-b_4x_3, b_4x_5-b_5x_1,b_5x_3-b_3x_5, b_6x_4-b_5x_6)$ and
$$\widetilde{\psi}^{T}=\left(\begin{array}{cccccc}
0 & b_4 & 0 & -b_1& 0 & 0\\
-b_2 & 0 & b_1 & 0 &0 & 0\\
0 & -b_3 & 0 & b_2 & 0 & 0\\
b_{3} & 0 & -b_4 &  0 & 0 & 0\\
-b_5&0&0&0&b_4&0\\
0&0&b_5&0&-b_3&0\\
0&0&0&b_6&0&-b_5
\end{array}\right) \left(\begin{array}{c} x_1\\
x_2\\
x_3\\
x_4\\
x_5\\
x_6
\end{array}\right)=B {\underline{x}}^T.$$
One can show that  $I_6(B)=b_5(b_1b_3-b_2b_4)^2(b_3,b_4,b_5)$.
In particular, if $b_5=0$ then $I_6(B)=0$ and thus no maximal submatrix of
$B$ is invertible. 
 
 Notice that  $b_1b_3 \neq
  b_2b_4$ by Corollary~\ref{productofbs}.
  Therefore $b_3$ and $b_4$ can not be
 simultaneously zero. Thus when $b_5
 \neq 0$ then $I_6(B) \neq 0$ and therefore $B$ has an invertible $6
 \times 6$ submatrix and
 $J':I=\m$. Hence $J \subset J'+\m I=J'$ and thus $J=J'$ and $J:I=\m$.

Suppose that $J$ is a general  minimal reduction of $I$, i.e. $J$ is generated by $\ell$ general elements of $I$. Then $a_j \in \mathbb{Q}$ and thus $b_j=a_j$ for all $ 1 \leq j \leq 6$. When $J$ is a  general minimal reduction we may choose $b_5 \neq 0$  and thus $J:I=\m$ for all such  $J$. Hence $\m I \subset J$ for all general minimal reductions $J$ of $I$. Therefore 
$\m I \subset \bigcap \limits_{J \in \mathcal{M}(I)} J$,
 where $ \mathcal{M}(I)= \{ J \mid J \mbox{ general minimal reduction of } I \}$. But we already saw that $\m I \not \subset \core{I}$ and therefore 
$\core{I} \neq \bigcap \limits_{J \in \mathcal{M}(I)} J$.

Finally, it is straightforward to see that $I$ satisfies $G_{\ell}$. If $I$ were weakly $(\ell-1)$-residually $S_2$ then by \cite[Theorem~4.5]{CPU01} the core would have been a finite intersection of general minimal reductions, a contradiction. Thus $I$ is not weakly $(\ell-1)$-residually $S_2$. Note that this can also be verified directly. \qed

\begin{remark} {\rm
Notice  that   Example~\ref{counterexample} establishes
that the condition that $I$ is weakly $(\ell-1)$-residually $S_2$ is necessary in  \cite[Theorem~4.5]{CPU01}.}
\end{remark}

\section{The core as a finite intersection}\label{finiteintersection}

We conclude this article by revisiting the question of whether the
core may be obtained as a finite intersection of minimal
reductions. Recall that under suitable assumptions Corso, Polini and
Ulrich prove that the core may be obtained as a finite intersection of
general minimal reductions \cite[Theorem~4.5]{CPU01}. Note that
Example~\ref{counterexample} is an instance where the assumptions
of \cite[Theorem~4.5]{CPU01} fail to hold and the core is not a
(finite) intersection of general minimal reductions.
We will prove in this
section that when $I$ is the edge ideal of an even cycle, then
$\core{I}$ is obtained as a finite intersection of minimal reductions
and we will give an explicit description of these minimal reductions. 
We first show that the edge ideal corresponding to an octagon is not weakly $(\ell-1)$-residually $S_2$.

\begin{example} \label{octagon}
Let   $I$ be the edge ideal of an even cycle of length $8$. Let $R$ be the corresponding localized polynomial ring over $\mathbb{Q}$. Then $I$ is not weakly $(\ell-1)$-residually $S_2$.\end{example}

\proof
Let $I=(e_1, \ldots, e_8)$. Then $\ell=7$. Let $\mathfrak{a}=(e_1+e_7-e_8, e_2+e_7+3e_8, e_3+e_7+e_8, e_4+e_7+e_8,e_5+e_7+e_8,e_6+e_7+2e_8)$ and $K=\mathfrak{a}:I$. Then ${\rm ht} \; K=6$ and ${\rm ht} (I+K)=7$. Therefore $K$ is a geometric $6$-residual intersection of $I$. Using Macaulay~2 \cite{M2} we have that ${\rm projdim}(R/K)=7$ and thus $\depth R/K=1$, which then means $R/K$ does not satisfy Serre's condition $S_2$. \QED

When  $I$ is the edge ideal of an even cycle then $I$ need  not be weakly $(\ell-1)$-residually $S_2$ as Example~\ref{octagon} suggests. Thus we may not apply \cite[Theorem~4.5]{CPU01}. Instead, we will employ different methods.

\begin{notation}\label{red for finite inter}{\rm
Let $I=(e_1, \ldots, e_d)$ be the edge ideal of an even cycle. 
For every $1 \leq t \leq d$, let $J_t =(e_1+a_1e_t, e_2+a_2e_t, \ldots , e_d+a_de_t)$, where   $a_i=1$ for $i \not=t$ and $a_t=-1$. For every $1 \leq t \leq d/2$ we define the following ideals:

$L_{2t}=(e_1+a_1e_{2t}, \ldots, e_d+a_de_{2t})$, where $a_i=1$ for all $i \neq 2t$ even, $a_i=0$ for all $i$ odd, and $a_{2t}=-1$;

$H_{2t}=(e_1+a_1e_{2t},  \ldots, e_d+a_de_{2t})$, where $a_i=1$ for all $i$ odd, $a_i=0$ for all $i \neq 2t$ even, and $a_{2t}=-1$;

$H_{2t-1}=(e_1+a_1e_{2t-1}, \ldots, e_d+a_de_{2t-1})$, where $a_i=1$ for all $i$ even, $a_i=0$ for all $i \neq 2t-1$ odd, and $a_{2t-1}=-1$. }

\end{notation}

\begin{remark} {\rm
Let $I$ be the edge ideal of an even cycle $e_1, \ldots ,e_d$. Using the same techniques as
in Examples~\ref{oddreduction} and~\ref{evenreduction}, we see that for every $1 \leq t \leq d/2$, the ideals $L_{2t}$, $H_{2t}$, $H_{2t-1}$ in~\ref{red for finite inter} are minimal reductions of $I$. When ${\rm char} \; k \neq 2$ then $J_t$ is a minimal reduction of $I$ for every $1\leq t \leq d$, by Example~\ref{basicreduction}. 
}
\end{remark}
 
\begin{proposition}\label{char not 2}
Let $I$ be the edge ideal of an even cycle $e_1, \ldots ,e_d$. Let $d=2n$ for some integer $n \geq 2$. Let $p={\rm char} \; k  \geq 0$. If $p \neq 2$ and  $n \not \equiv 1 \mod p$  then $\core{I}=\bigcap\limits_{t=1}^d J_t$.
\end{proposition}

\proof
First recall that $\core{I}=\m I$, by Theorem~\ref{formulaext}. Let $C=\bigcap \limits_{t=1}^d
J_t$. Since
$J_t$ is a minimal reduction of $I$ for each $t$, we have that
$\m I\subset C$.

In order to establish the other inclusion,  suppose $f\in C 
\setminus \m I$. Since $f \in I$ then we may write $f=\sum \limits_{i=1}^d h_ie_i$, for some $h_i
  \in R$. Now since (by clearing denominators in the localization if
  necessary) $h_i$ can be taken to be a polynomial then we may write
  $h_i=g_i+h_i'$, for some  $g_i\in \m$ and $h_i' \in k$ of degree
  $0$. Notice that $g_ie_i \in 
  \m I \subset J_t$ for all $i,t$. Thus if $g=\sum \limits_{i=1}^d
  g_ie_i$, then $g \in \m I \subset C$ and so
  $f'=f-g=\sum \limits_{i=1}^d h_i'e_i \in C \setminus \m I$. Therefore, without
  loss of generality, we may assume $f=\sum \limits_{i=1}^d h_ie_i$, where
  $h_1=1$ and $h_i \in k$ for all $i$.

We observe that since  $f\in J_1$ then we may write $f=\sum \limits_{i=2}^{d} a_i
(e_i+e_1)$, for some $a_i \in R$. Notice that  $f$ is homogeneous of degree
 $2$ and thus we may assume $a_i 
\in k$ since all terms of higher degree must cancel. The set $\{ e_1, \ldots, e_d \}$
is linearly independent over $k$. Therefore we may equate coefficients
of $e_i$ in the two summation representations of $f$. Thus $h_i=a_i$
for $i \geq 2$. Furthermore, by equating the coefficients of $e_1$ we have
$\sum \limits_{i=2}^d h_i =1$.
Since $f\in
J_2$, then $f=b_1(e_1+e_2) + \sum \limits_{i=3}^d b_i(e_i+e_2)$, for some
$b_i \in k$. Using the same method as above we obtain 
$b_i=h_i$ for all $i \not= 2$.  By examining $e_2$, and recalling that
$h_1=1$, we see that $1 + 
\sum \limits_{i=3}^d h_i = h_2$ and thus  $1+ \sum \limits_{i=2}^d h_i =2h_2$. Combining both equations yields $2=2h_2$. Since $p \neq 2$, we have $h_2 = 1$ and $\sum \limits_{i=3}^d h_i = 0$.

We will proceed by induction.  Suppose that for some $t<d$, $h_i=1$ for  all $i \leq t$ and
$\sum \limits_{i=t+1}^d h_i \equiv 2-t$, where equivalence will be considered modulo $p$. Since $f \in J_{t+1}$ then $f= \sum \limits_{i=1}^{t}(e_i+ e_{t+1})+ \sum \limits_{i=t+2}^d h_i(e_i+e_{t+1})$. Examining
    the coefficient of $e_{t+1}$ yields 
    $h_{t+1} \equiv t+\sum \limits_{i=t+2}^d h_i $ and therefore
    $2h_{t+1}\equiv t+\sum \limits_{i=t+1}^d h_i $, or $
    2h_{t+1}\equiv t+2-t$. So $h_{t+1}\equiv 1$
    and thus $\sum \limits_{i=t+2}^d h_i \equiv 1-t=2-(t+1)$. Thus by induction, we may assume $h_i \equiv 1$ for all $i \leq d-1$ and
$\sum \limits_{i=t+1}^d h_i  \equiv 2-t $ for all $t \leq
    d-1$. Note that since $h_i \in k$, $h_i \equiv 1$ implies
    $h_i=1$ in $k$. 
Now assume $t=d-1$. Then 
$h_d=\sum \limits_{i=d}^d h_i \equiv 2-(d-1) = 3-d $. Again, since  $f\in J_d$ then 
$f=\sum \limits_{i=1}^{d-1}(e_i+e_d)$ and thus $h_d=d-1$. 
But $h_d \equiv 3-d $, so $d-1\equiv 3-d $ or $d\equiv 2$. Equivalently, since $d=2n$ then $n \equiv 1$, which is a contradiction. Therefore $C \subset \m I$. \QED

We now consider the remaining  cases when the characteristic of the residue field is $2$ or $n  \equiv 1 \mod p$.

\begin{proposition} \label{char 2}
Let $I$ be the edge ideal of an even cycle $e_1, \ldots, e_d$. Let $d=2n$ for some integer $n \geq 2$. Let $p={\rm char} \; k  \geq 0$. If  $p=2$ or $n  \equiv 1 \mod p$ then  $$\core{I}=\bigcap\limits_{i=1}^d H_i  \cap \bigcap \limits_{t=1}^{d/2} L_{2t},$$ where $H_{i}$ and $L_{2t}$ are as in \ref{red for finite inter}. 
\end{proposition}

\proof
Let $C=\bigcap\limits_{i=1}^d H_i  \cap \bigcap \limits_{t=1}^{d/2} L_{2t}$.
By Theorem~\ref{formulaext}, $\core{I}=\m I$. Since for every $ 1 \leq t \leq d/2$ and every $1 \leq i \leq d$ we have that $L_{2t}$ and $H_i$ are all minimal reductions of $I$, then $\m I\subset C $. 
As before, we may assume $f \in C \setminus \m I$ and $f=\sum
  \limits_{i=1}^d h_ie_i$, where 
  $h_1=1$ and $h_i \in k$ for all $i$. 

First we note that since  $f\in H_1$  we may write $f= \sum
\limits_{i=1}^{d/2}a_{2i}(e_{2i}+e_1)+ \sum
\limits_{i=2}^{d/2}a_{2i-1} e_{2i-1}$ for some $a_i \in k$. Equating coefficients yields
$a_i=h_i$ for all $i\not= 1$ and that $\sum \limits_{i=1}^{d/2}a_{2i}=\sum \limits_{i=1}^{d/2}h_{2i}=1$.

Since $f \in L_d$ then $f= \sum
\limits_{i=1}^{d/2-1}b_{2i}(e_{2i}+e_d)+ \sum
\limits_{i=1}^{d/2}b_{2i-1} e_{2i-1}$ for some $b_i \in k$. Equating
coefficients as before, 
we have  that  $\sum \limits_{i=1}^{d/2-1}
h_{2i}=h_d$ and thus $\sum \limits_{i=1}^{d/2} h_{2i}=2h_d$. Hence
$2h_d=1$. If ${\rm char} \;k=2$ then we have that $0=1$, which is a
contradiction. Thus we may assume that ${\rm char} \; k \neq 2$, $n
\equiv 1 \mod p$ and $2h_d=1$. 

Similarly, since $f\in L_{d-2}$ we obtain $\sum \limits_{i=1}^{d/2-2} h_{2i}+h_d=h_{d-2}$ and hence $\sum \limits_{i=1}^{d/2} h_{2i}=2h_{d-2}$. Thus $2h_{d-2}=1$. Repeating this process yields $2h_{2i}=1$ for all $1 \leq i \leq d/2$. But as $\sum \limits_{i=1}^{d/2} h_{2i}=1$ we have $2\sum \limits_{i=1}^{d/2} h_{2i}= {\frac{d}{2}}=2$, i.e. $d \equiv4 \mod p$. Since $d=2n$  then $n \equiv 2 \mod p$, which is a contradiction. Thus $C \subset \m I$.
\qed

\begin{theorem} \label{finite intersection core}
Let $I$ be the edge ideal of an even cycle. Then $\core{I}$ is obtained as a finite intersection of minimal reductions of $I$.
\end{theorem}

\proof Combine Propositions~\ref{char not 2} and~\ref{char 2}. \qed

\begin{remark} \label{gradecore}{\rm Let $I$ be the edge ideal of an even cycle. Recall
    that the ${\rm gradedcore} (I)$ is the intersection of all homogeneous
    minimal reductions of $I$.  In general, $\core{I} \subset {\rm
      gradedcore}(I) $. 
     We  note that all the reductions in \ref{red for finite inter} are homogeneous minimal reductions. Hence ${\rm gradedcore}(I) \subset C$,  
      where $C$ is as in
    Propositions~\ref{char not 2} and \ref{char 2}. Therefore,
    $\core{I}={\rm gradedcore}(I)=\m I$, by Theorem~\ref{finite
      intersection core}.}
\end{remark}

\end{document}